\documentclass[a4paper, 11pt]{article}

\usepackage{fullpage}

\usepackage[english]{babel}
\usepackage[utf8]{inputenc}
\usepackage{amsmath, amsfonts, amssymb}
\usepackage{mathtools}
\usepackage{bm}
\usepackage{authblk}

\usepackage{booktabs}
\usepackage{graphicx}
\usepackage{subcaption}

\usepackage[linesnumbered,ruled,vlined]{algorithm2e}
\SetKwInput{KwInput}{input}
\SetAlFnt{\small}

\newcommand{\Real}{\mathbb{R}}

\title{Note: low-rank tensor train completion with side information based on Riemannian optimization}
\author[1,2]{Stanislav Budzinskiy}
\author[1]{Nikolai Zamarashkin}
\affil[1]{Marchuk Institute of Numerical Mathematics of the Russian Academy of Sciences}
\affil[2]{Faculty of Computational Mathematics and Cybernetics, Lomonosov Moscow State University}
\date{}

\begin{document}
\maketitle

\begin{abstract}
We consider the low-rank tensor train completion problem when additional side information is available in the form of subspaces that contain the mode-$k$ fiber spans. We propose an algorithm based on Riemannian optimization to solve the problem. Numerical experiments show that the proposed algorithm requires far fewer known entries to recover the tensor compared to standard tensor train completion methods.
\end{abstract}

\section{Introduction}
Let $\bm{\mathcal{A}} \in \Real^{N_1 \times \dots \times N_d}$ be a $d$-dimensional tensor with low tensor train (TT) ranks \cite{OseledetsTensorTrain2011}. We consider the problem of recovering $\bm{\mathcal{A}}$ from a limited number of its entries. Put more formally, let
\[ 
    \Omega \subset [N_1] \times \ldots \times [N_d]
\]
be the indices of the given entries and let
\[
    P_\Omega : \Real^{N_1 \times \dots \times N_d} \to \Real^{N_1 \times \dots \times N_d}
\]
be the projection operator defined by
\[
    P_\Omega \bm{\mathcal{A}}(i_1,\ldots,i_d) = 
    \begin{cases}
        \bm{\mathcal{A}}(i_1,\ldots,i_d), & (i_1,\ldots,i_d) \in \Omega,\\
        0, & (i_1,\ldots,i_d) \not\in \Omega.
    \end{cases}
\]
If the TT ranks of $\bm{\mathcal{A}}$ are equal to $\bm{r} = (1, r_1, \ldots, r_{d - 1}, 1)$, the TT completion problem can be written down as
\begin{equation}\label{eq:completion_no_side}
    f(\bm{\mathcal{X}}) = \frac{1}{2} \| P_\Omega \bm{\mathcal{X}} - P_\Omega \bm{\mathcal{A}} \|^2_F \to \min_{\bm{\mathcal{X}}} \quad\text{s.t.}\quad \mathrm{rank}_{TT}(\bm{\mathcal{X}}) = \bm{r}.
\end{equation}
It is known that the set
\[
    \mathcal{N}_{\bm{r}} = \{ \bm{\mathcal{X}} \in \Real^{N_1 \times \dots \times N_d} : \mathrm{rank}_{TT}(\bm{\mathcal{X}}) = \bm{r} \}
\]
is a smooth embedded submanifold of $\Real^{N_1 \times \dots \times N_d}$, hence the completion problem can be solved via Riemannian optimization \cite{SteinlechnerRiemannian2016}. Other approaches include alternating least squares \cite{GrasedyckEtAlVariants2015, GrasedyckKramerStable2019} and nuclear norm minimization \cite{BenguaEtAlEfficient2017}.

Sometimes, additional side information about $\bm{\mathcal{A}}$ is available that could be incorporated into the algorithm to reduce the amount of data $|\Omega|$ needed for successful reconstruction. This side information can come in the form of orthogonal matrices
\[
    \bm{Q}_k \in \Real^{N_k \times M_k}, \quad \bm{Q}_k^T \bm{Q}_k = \bm{I}_{M_k}, \quad k = 1,\ldots,d,
\]
that constraint the mode-$k$ subspaces of $\bm{\mathcal{A}}$ as
\begin{equation}\label{eq:side_inf_colspan}
    \mathrm{col}(\bm{A}_{(k)}) \subset \mathrm{col}(\bm{Q}_k), \quad k = 1,\ldots,d,
\end{equation}
where $\bm{A}_{(k)} \in \Real^{N_k \times (N_1\ldots N_{k-1}N_{k+1}\ldots N_d)}$ is the mode-$k$ flattening of $\bm{\mathcal{A}}$ and $\mathrm{col}(\cdot)$ is the column span of the matrix argument. Put differently,
\[
    \bm{Q}_k \bm{Q}_k^T \bm{A}_{(k)} = \bm{A}_{(k)}, \quad k = 1,\ldots,d,
\]
or
\[
    \bm{\mathcal{A}} \times_1 (\bm{Q}_1 \bm{Q}_1^T) \times_2 \ldots \times_{d} (\bm{Q}_d \bm{Q}_d^T) = \bm{\mathcal{A}},
\]
where $\times_k$ is the $k$-mode product. 

To the best of our knowledge, low-rank completion problems with side information have been analyzed only for the matrix case \cite{JainDhillonProvable2013, XuEtAlSpeedup2013, ChiangEtAlMatrix2015}; the corresponding algorithms have been applied in bioinformatics \cite{NatarajanDhillonInductive2014}. For the nuclear norm minimization algorithm in \cite{XuEtAlSpeedup2013}, it was proved that it requires only $O(r M \log M \log N)$ known entries to recover the matrix instead of $O(r N \log^2 N)$ required by the standard nuclear norm minimization \cite{RechtSimpler2011}. Note that $O(r N \log^2 N)$ is the bound for a particular algorithm; meanwhile, $O(r N \log N)$ samples are necessary for any matrix completion method whatsoever \cite{CandesTaoPower2010}.

The goal of this paper is to develop a Riemannian optimization  algorithm to solve the low-rank TT completion problem \eqref{eq:completion_no_side} given the side information \eqref{eq:side_inf_colspan} and to examine how the side information affects the number of known entries needed for successful recovery.

\section{Structure of the feasible set}
Let 
\[
    \mathcal{N}_{\bm{r}}(\bm{Q}) = \{ \bm{\mathcal{X}} \in \mathcal{N}_{\bm{r}} : \bm{Q}_k \bm{Q}_k^T \bm{X}_{(k)} = \bm{X}_{(k)},\,k = 1,\ldots,d \}
\]
be the set of low-rank TTs that conform to the side information. Each $\bm{\mathcal{X}} \in \mathcal{N}_{\bm{r}}(\bm{Q})$ can be represented as
\[
    \bm{\mathcal{X}} = Q \bm{\mathcal{Y}} \equiv \bm{\mathcal{Y}} \times_1 \bm{Q}_1 \times_2 \ldots \times_d \bm{Q}_d
\]
where $\bm{\mathcal{Y}} \in \Real^{M_1 \times \ldots \times M_d}$ is a smaller tensor with the same TT ranks:
\[
    \bm{\mathcal{Y}} \in \mathcal{M}_{\bm{r}} = \{ \bm{\mathcal{Y}} \in \Real^{M_1 \times \dots \times M_d} : \mathrm{rank}_{TT}(\bm{\mathcal{Y}}) = \bm{r} \}.
\]
Owing to the orthogonality of $\{ \bm{Q}_k \}$, the mapping
\[
    Q: \Real^{M_1 \times \dots \times M_d} \to \Real^{N_1 \times \dots \times N_d}
\]
establishes a one-to-one correspondence between $\mathcal{M}_{\bm{r}}$ and $\mathcal{N}_{\bm{r}}(\bm{Q})$. Its inverse acts according to
\[
    Q^{-1} \bm{\mathcal{X}} = Q^{T} \bm{\mathcal{X}} \equiv \bm{\mathcal{X}} \times_1 \bm{Q}_1^T \times_2 \ldots \times_d \bm{Q}_d^T.
\]

\section{Two ways to use the one-to-one correspondence}
This correspondence lies at the heart of the methods in the matrix case \cite{JainDhillonProvable2013, XuEtAlSpeedup2013}. With its help, we can formulate the TT completion problem with side information as
\begin{equation}\label{eq:completion_side_bad}
    g(\bm{\mathcal{Y}}) = \frac{1}{2} \| P_\Omega Q \bm{\mathcal{Y}} - P_\Omega \bm{\mathcal{A}} \|^2_F \to \min_{\bm{\mathcal{Y}}} \quad\text{s.t.}\quad \bm{\mathcal{Y}} \in \mathcal{M}_{\bm{r}}.
\end{equation}
Since the set $\mathcal{M}_{\bm{r}}$ is a smooth embedded submanifold of $\Real^{M_1 \times \ldots \times M_d}$, the optimization problem can be solved with Riemannian gradient descent or Riemannian conjugate gradient methods \cite{SteinlechnerRiemannian2016}. At each iteration, Riemannian optimization methods perform a series of steps:
\begin{enumerate}
    \item compute the usual Euclidean gradient $\nabla g(\bm{\mathcal{Y}})$,
    \item compute the Riemannian gradient $\mathrm{grad}~ g(\bm{\mathcal{Y}})$ by projecting the Euclidean gradient $\nabla g(\bm{\mathcal{Y}})$ onto the tangent space $T_{\bm{\mathcal{Y}}} \mathcal{M}_{\bm{r}}$,
    \item choose the step direction $\bm{\eta} \in T_{\bm{\mathcal{Y}}} \mathcal{M}_{\bm{r}}$ and the step size $\alpha \in \Real$,
    \item project $\bm{\mathcal{Y}} + \alpha\bm{\eta}$ onto the manifold with a retraction.
\end{enumerate}
A computationally efficient way to evaluate the Riemannian gradient for the low-rank TT completion problem was proposed in \cite{SteinlechnerRiemannian2016}. It exploits the sparsity of the Euclidean gradient and takes $O(d|\Omega|r^2 + dMr^3)$ operations. However, the Euclidean gradient of \eqref{eq:completion_side_bad} is no longer sparse due to the presence of the mapping $Q$. Namely,
\[
    \nabla g(\bm{\mathcal{Y}}) = (P_\Omega Q \bm{\mathcal{Y}} - P_\Omega \bm{\mathcal{A}}) \times_1 \bm{Q}_1^T \times_2 \ldots \times_d \bm{Q}_d^T.
\]
While the tensor $P_\Omega Q \bm{\mathcal{Y}} - P_\Omega \bm{\mathcal{A}}$ is indeed sparse, the sparsity is broken by the mode-products and the projection onto the tangent space thus suffers from the curse of dimensionality and becomes unacceptably inefficient. 

This is why we suggest to use the correspondence between $\mathcal{M}_{\bm{r}}$ and $\mathcal{N}_{\bm{r}}(\bm{Q})$ in a different manner. Since $\mathcal{M}_{\bm{r}}$ is a smooth submanifold of $\Real^{M_1 \times \ldots \times M_d}$ and $Q$ is a one-to-one linear mapping, $\mathcal{N}_{\bm{r}}(\bm{Q}) = Q\mathcal{M}_{\bm{r}}$ is also a smooth submanifold but in the larger space $\Real^{N_1 \times \ldots \times N_d}$. This means that we can pose the low-rank TT completion problem with side information as
\begin{equation}\label{eq:completion_side_good}
    f(\bm{\mathcal{X}}) = \frac{1}{2} \| P_\Omega \bm{\mathcal{X}} - P_\Omega \bm{\mathcal{A}} \|^2_F \to \min_{\bm{\mathcal{X}}} \quad\text{s.t.}\quad \bm{\mathcal{X}} \in \mathcal{N}_{\bm{r}}(\bm{Q}).
\end{equation}

\section{Riemannian optimization scheme}
To apply Riemannian optimization methods to problem \eqref{eq:completion_side_good}, we need to describe the tangent space $T_{\bm{\mathcal{X}}} \mathcal{N}_{\bm{r}}(\bm{Q})$, the projection operator $P_{T_{\bm{\mathcal{X}}} \mathcal{N}_{\bm{r}}(\bm{Q})}: \Real^{N_1 \times \ldots \times N_d} \to T_{\bm{\mathcal{X}}} \mathcal{N}_{\bm{r}}(\bm{Q})$ onto the tangent space, and the retraction operator $R: T \mathcal{N}_{\bm{r}}(\bm{Q}) \to \mathcal{N}_{\bm{r}}(\bm{Q})$ from the tangent bundle  into the submanifold. In the case of the Riemannian conjugate gradient method, the vector transport
\[
    \mathcal{T}_{\bm{\mathcal{X}} \to \bm{\mathcal{X}}_+}: T_{\bm{\mathcal{X}}} \mathcal{N}_{\bm{r}}(\bm{Q}) \to T_{\bm{\mathcal{X}}_+} \mathcal{N}_{\bm{r}}(\bm{Q})
\]
needs to be provided too.

To describe $T_{\bm{\mathcal{X}}} \mathcal{N}_{\bm{r}}(\bm{Q})$ we will use the correspondence between $\mathcal{N}_{\bm{r}}(\bm{Q})$ and $\mathcal{M}_{\bm{r}}$. Since the mapping $Q$ is linear, the tangent space is
\[
    T_{\bm{\mathcal{X}}} \mathcal{N}_{\bm{r}}(\bm{Q}) = Q T_{Q^{-1}\bm{\mathcal{X}}} \mathcal{M}_{\bm{r}}.
\]
We can further prove that
\[
    T_{\bm{\mathcal{X}}} \mathcal{N}_{\bm{r}}(\bm{Q}) = Q Q^T T_{\bm{\mathcal{X}}} \mathcal{N}_{\bm{r}}.
\]
This immediately gives the projection operator onto $T_{\bm{\mathcal{X}}} \mathcal{N}_{\bm{r}}(\bm{Q})$ as
\begin{equation}\label{eq:proj_tan}
    P_{T_{\bm{\mathcal{X}}} \mathcal{N}_{\bm{r}}(\bm{Q})} = Q Q^T P_{T_{\bm{\mathcal{X}}} \mathcal{N}_{\bm{r}}},
\end{equation}
and the description of $P_{T_{\bm{\mathcal{X}}} \mathcal{N}_{\bm{r}}}$ together with an efficient way to apply it are presented in \cite{SteinlechnerRiemannian2016}.

In the absence of side information, TT-SVD \cite{OseledetsTensorTrain2011} defines a retraction \cite{SteinlechnerRiemannian2016} as
\[
     R : T \mathcal{N}_{\bm{r}} \ni (\bm{\mathcal{X}}, \bm{\eta}) \mapsto P_{\bm{r}}^{TT}(\bm{\mathcal{X}} + \bm{\eta}) \in \mathcal{N}_{\bm{r}},
\]
where $P_{\bm{r}}^{TT}$ maps a tensor to its rank-$\bm{r}$ TT approximation. The TT-SVD algorithm and its TT-rounding variant perform sequences of QR decompositions and truncated SVDs that cannot enlarge the fiber spans, i.e. if $\mathrm{col}(\bm{\mathcal{X}}_{(k)}) \subset \mathrm{col}(\bm{Q}_k)$ then $\mathrm{col}((P_{\bm{r}}^{TT}\bm{\mathcal{X}})_{(k)}) \subset \mathrm{col}(\bm{Q}_k)$. This means that if $\bm{\mathcal{X}} \in \mathcal{N}_{\bm{r}}(\bm{Q})$ and $\bm{\eta} \in T_{\bm{\mathcal{X}}} \mathcal{N}_{\bm{r}}(\bm{Q})$, then $\mathrm{col}((\bm{\mathcal{X}} + \bm{\eta})_{(k)}) \subset \mathrm{col}(\bm{Q}_k)$ and $P_{\bm{r}}^{TT}(\bm{\mathcal{X}} + \bm{\eta}) \in \mathcal{N}_{\bm{r}}(\bm{Q})$. Thus the function 
\[
    R(\bm{Q}): T \mathcal{N}_{\bm{r}}(\bm{Q}) \ni (\bm{\mathcal{X}}, \bm{\eta}) \mapsto P_{\bm{r}}^{TT}(\bm{\mathcal{X}} + \bm{\eta}) \in \mathcal{N}_{\bm{r}}(\bm{Q})
\]
defines a retraction since it inherits the properties of $R$.

As in \cite{SteinlechnerRiemannian2016}, the orthogonal projection $P_{T_{\bm{\mathcal{X}}_+} \mathcal{N}_{\bm{r}}}(\bm{Q})$ can be used as the vector transport for the nonlinear conjugate gradient scheme. But since, just like TT-SVD, $P_{T_{\bm{\mathcal{X}}_+} \mathcal{N}_{\bm{r}}}$ does not enlarge the fiber spans and tensors in $T_{\bm{\mathcal{X}}} \mathcal{N}_{\bm{r}}(\bm{Q})$ conform to the side information, it suffices to apply $P_{T_{\bm{\mathcal{X}}_+} \mathcal{N}_{\bm{r}}}$.

This being said, we can solve the low-rank TT completion problem with side information \eqref{eq:completion_side_good} with the RTTC algorithm of \cite{SteinlechnerRiemannian2016} by changing the projection operator onto the tangent space as in \eqref{eq:proj_tan} (provided that the initial approximation for the iterations has correct fiber spans). We call this algorithm RTTC with side information (RTTC-SI).

\section{Numerical experiments}
We test the algorithm on random synthetic data. At first, a tensor $\bm{\mathcal{\tilde{A}}} \in \Real^{N \times \ldots \times N}$ of order $d$ is generated as a TT of rank $\bm{r} = (1, r, \ldots, r, 1)$ with standard Gaussian TT cores. For each mode, we generate a random standard Gaussian matrix $\bm{\tilde{Q}}_k \in \Real^{N \times M}$, $k = 1,\ldots,d$, and orthogonalize it via Gram-Schmidt to get $\bm{Q}_k \in \Real^{N \times M}$, $\bm{Q}_k^T \bm{Q}_k = \bm{I}_M$. We then form the tensor to be recovered as
\[
    \bm{\mathcal{A}} = Q Q^T \bm{\mathcal{\tilde{A}}} \in \Real^{N \times \ldots \times N}.
\]
The set of indices $\Omega$ for which the entries are known is generated uniformly at random. The initial approximation $\bm{\mathcal{X}}_0 \in \Real^{N \times \ldots \times N}$ is generated at random similarly to $\bm{\mathcal{A}}$:
\[
    \bm{\mathcal{X}}_0 = Q Q^T \bm{\mathcal{\tilde{X}}}_0 \in \Real^{N \times \ldots \times N}
\]
with a standard Gaussian rank-$\bm{r}$ TT $\bm{\mathcal{\tilde{X}}}_0$.

We let the algorithm take 250 iterations and call the recovery convergent if the relative error is small on a test set:
\[
    \frac{\| P_\Gamma \bm{\mathcal{X}}_{250} - P_\Gamma \bm{\mathcal{A}} \|_F}{\| P_\Gamma \bm{\mathcal{A}} \|_F} < 10^{-4}
\]
for $\Gamma \in [N]^d$, $|\Gamma| = |\Omega|$, generated uniformly at random. For each combination of parameters we run 5 experiments and report the frequency of convergent runs.

\begin{figure}[h]
\begin{subfigure}[b]{0.5\textwidth}
\centering
	\includegraphics[width=1\textwidth]{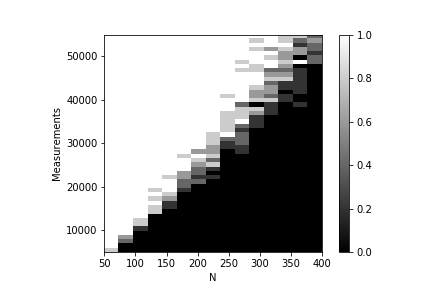}
    \caption{RTTC}
\end{subfigure}%
\begin{subfigure}[b]{0.5\textwidth}
\centering
	\includegraphics[width=1\textwidth]{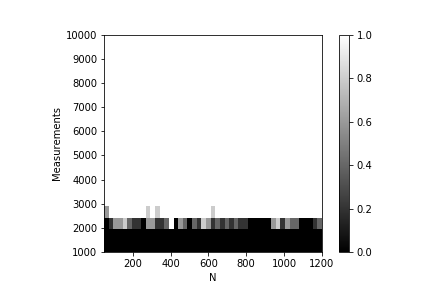}
    \caption{RTTC-SI}
\end{subfigure}
\caption{Phase plots of RTTC (a) and RTTC-SI (b) for $d = 5$, $r = 3$, $M = 30$, and varying $N$.}
\label{fig:cmplt_m30}
\end{figure}

In Figure \ref{fig:cmplt_m30}, we compare the phase plots obtained with RTTC and RTTC-SI. The phase plot of RTTC is identical to that of \cite{SteinlechnerRiemannian2016} which shows that the algorithm is ignorant to the side information even when the initial approximation has the correct fiber spans. At the same time, RTTC-SI converges even when the available data are scarce and the sufficient number of entries does not depend on $N$.

\begin{figure}[h]
\begin{subfigure}[b]{0.5\textwidth}
\centering
	\includegraphics[width=1\textwidth]{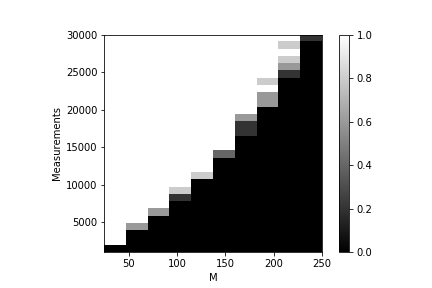}
	\caption{$N = 250$}
\end{subfigure}%
\begin{subfigure}[b]{0.5\textwidth}
\centering
	\includegraphics[width=1\textwidth]{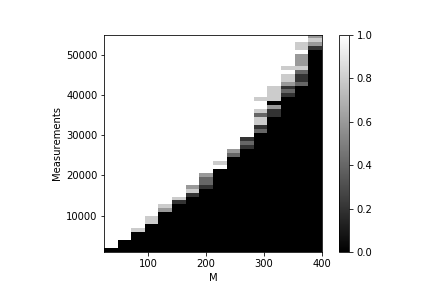}
    \caption{$N = 400$}
\end{subfigure}
\caption{Phase plots of RTTC-SI for $d = 5$, $r = 3$, and varying $M$.}
\label{fig:cmplt_n}
\end{figure}

\begin{figure}[h]
\begin{subfigure}[b]{0.5\textwidth}
\centering
    \includegraphics[width=\textwidth]{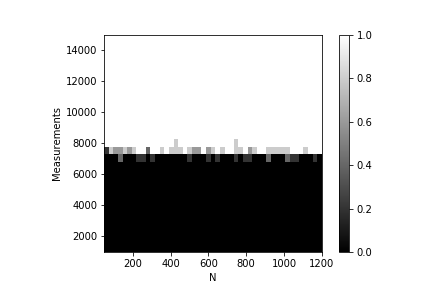}
    \caption{$d = 5$, $r = 6$}
\end{subfigure}%
\begin{subfigure}[b]{0.5\textwidth}
\centering
    \includegraphics[width=\textwidth]{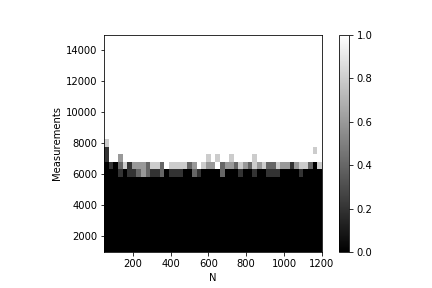}
    \caption{$d = 10$, $r = 3$}
\end{subfigure}%
\caption{Phase plots of RTTC-SI for $M = 30$, and varying $N$.}
\label{fig:cmplt_rd}
\end{figure}

The phase plots in Figure \ref{fig:cmplt_n} illustrate that phase transition happens when $|\Omega| = O(M \log M)$ similarly to the matrix case. Figure \ref{fig:cmplt_rd} allows to see how the phase transition curve changes as we double $r$ or $d$ (compare with Figure \ref{fig:cmplt_m30} (b)).

\begin{figure}
\centering
    \includegraphics[width=0.5\textwidth]{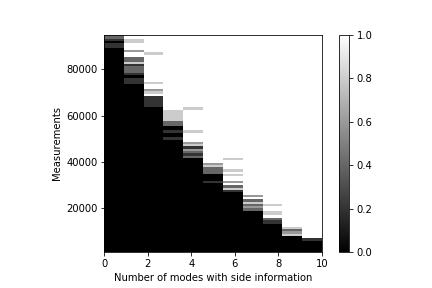}
    \caption{Phase plots of RTTC-SI for $N = 200$, $d = 10$, $r = 3$, $M = 30$ when the side information is available for a limited number of modes.}
    \label{fig:cmplt_m30d10_nmodes}
\end{figure}

Finally, we consider the situation when the side information is available for a limited number of modes. Let 
\[  
    M_k =
    \begin{cases}
        M, & k \leq l,\\
        N, & k > l,
    \end{cases} \quad k = 1,\ldots,d.
\]
In Figure \ref{fig:cmplt_m30d10_nmodes} we present the phase plots for $l$ varying from $0$ to $d$.

\section{Conclusion}
In this paper we presented the RTTC-SI algorithm for low-rank TT completion with side information and numerically studied the effects of side information on the phase transition curve.

\bibliography{main}
\bibliographystyle{unsrt}

\end{document}